\documentclass[11pt]{article}

\usepackage{amscd,amsmath, amssymb, fancyhdr, epsfig,url,MnSymbol}
\usepackage{color}
\usepackage[matrix,arrow]{xy}
\usepackage{diagrams}
\numberwithin{equation}{section}

\newcommand{\version}{version 2.0,\ \  April 8, 2016}

\def\eqref#1{(\ref{#1})}
\newcommand{\goth}{\mathfrak}
\newcommand{\g}{{\mathfrak g}}
\renewcommand{\a}{{\mathfrak a}}
\newcommand{\arrow}{{\:\longrightarrow\:}}
\newcommand{\Z}{{\Bbb Z}}
\newcommand{\C}{{\Bbb C}}

\newcommand{\R}{{\Bbb R}}
\newcommand{\Q}{{\Bbb Q}}

\def\1{\sqrt{-1}\:}

\makeatletter
\def\x@arrow{\DOTSB\Relbar}
\def\xlongequalsignfill@{\arrowfill@\x@arrow\Relbar\x@arrow}
\newcommand{\xlongequal}[2][]{%
        \ext@arrow 0099\xlongequalsignfill@{#1}{#2}}
\def\xlongrightarrowfill@{\arrowfill@\relbar\relbar\longrightarrow}
\newcommand{\xlongrightarrow}[2][]{%
        \ext@arrow 0099\xlongrightarrowfill@{#1}{#2}}
\makeatother

\renewcommand{\phi}{\varphi}
\renewcommand{\epsilon}{\varepsilon}
\renewcommand{\geq}{\geqslant}
\renewcommand{\leq}{\leqslant}
\renewcommand{\max}{{\rm max}}

\newcommand{\im}{\operatorname{im}}

\newcommand{\rank}{\operatorname{\sf rank}}
\newcommand{\Span}{\operatorname{\sf span}}
\newcommand{\Av}{\operatorname{\sf Av}}
\newcommand{\End}{\operatorname{End}}

\newcommand{\Alb}{\operatorname{Alb}}
\newcommand{\Mer}{\operatorname{Mer}}

\newcounter{Mycounter}[section]
\newcounter{lemma}[section]
\newcounter{claim}[section]
\renewcommand{\theclaim}{{Claim \thesection.\arabic{claim}}}
\newcommand{\claim}{%
    \setcounter{claim}{\value{Mycounter}}
    \refstepcounter{claim}
    \stepcounter{Mycounter}
    {\noindent \bf \theclaim:\ }}

\newcounter{sublemma}[section]
\newcounter{corollary}[section]
\renewcommand{\thecorollary}{{Corollary \thesection.\arabic{corollary}}}
\newcommand{\corollary}{%
    \setcounter{corollary}{\value{Mycounter}}
    \refstepcounter{corollary}
    \stepcounter{Mycounter}
    {\noindent \bf \thecorollary:\ }}

\newcounter{theorem}[section]
\renewcommand{\thetheorem}{{Theorem \thesection.\arabic{theorem}}}
\newcommand{\theorem}{%
    \setcounter{theorem}{\value{Mycounter}}
    \refstepcounter{theorem}
    \stepcounter{Mycounter}
    {\noindent \bf \thetheorem:\ }}

\newcounter{conjecture}[section]
\newcounter{proposition}[section]
\renewcommand{\theproposition}
      {{Proposition \thesection.\arabic{proposition}}}
\newcommand{\proposition}{%
    \setcounter{proposition}{\value{Mycounter}}
    \refstepcounter{proposition}
    \stepcounter{Mycounter}
    {\noindent \bf \theproposition:\ }}

\newcounter{definition}[section]
\renewcommand{\thedefinition}{{Definition~\thesection.\arabic{definition}}}
\newcommand{\definition}{%
    \setcounter{definition}{\value{Mycounter}}
    \refstepcounter{definition}
    \stepcounter{Mycounter}
    {\noindent \bf \thedefinition:\ }}

\newcounter{example}[section]
\renewcommand{\theexample}{{Example \thesection.\arabic{example}}}
\newcommand{\example}{%
    \setcounter{example}{\value{Mycounter}}
    \refstepcounter{example}
    \stepcounter{Mycounter}
    {\noindent \bf \theexample:\ }}

\newcounter{remark}[section]
\renewcommand{\theremark}{{Remark \thesection.\arabic{remark}}}
\newcommand{\remark}{%
    \setcounter{remark}{\value{Mycounter}}
    \refstepcounter{remark}
    \stepcounter{Mycounter}
    {\noindent \bf \theremark:\ }}

\newcounter{problem}[section]
\newcounter{question}[section]
\makeatletter

\@addtoreset{equation}{section} \@addtoreset{footnote}{section} \makeatother

\def\blacksquare{\hbox{\vrule width 5pt height 5pt depth 0pt}}
\def\endproof{\blacksquare}

\begin{document}
%%%%%%%%%%%%%%%%%%%%%%%%%%%%%%%%%%%%%%%%%%%%%%%%%%%%%%%%%%%%
\begin{center}
{\LARGE\bf
Algebraic dimension of complex nilmanifolds\\[4mm]
}
%%%%%%%%%%%%%%%%%%%%%%%%%%%%%%%%%%%%%%%%%%%%%%%%%%%%%%%%%%%%

Anna Fino, Gueo Grantcharov, Misha Verbitsky\footnote{Anna Fino
is partially supported by   PRIN, FIRB and GNSAGA (Indam), Gueo
Grantcharov is supported by a grant from the Simons
Foundation (\#246184), and Misha Verbitsky is
partially supported by RSCF grant 14-21-00053
within AG Laboratory NRU-HSE.}

\end{center}

%%%%%%%%%%%%%%%%%%%%%%%%%%%%%%%%%%%%%%%%%%%%%%%%

{\small
\begin{minipage}[t]{0.85\linewidth}
{\bf Abstract} \\
Let $M$ be a complex nilmanifold, that is, a compact quotient
of a nilpotent Lie group endowed with an invariant complex
structure  by a discrete lattice. A holomorphic
differential on $M$ is a closed, holomorphic 1-form. We show that $a(M)\leq k$,
where  $a(M)$ is the  algebraic dimension $a(M)$ 
(i.e.  the transcendence degree of the
field of meromorphic functions) and  $k$ is the dimension of the space
of holomorphic differentials. We prove a similar result about 
meromorphic maps to Kahler manifolds.
\end{minipage}
}

%%%%%%%%%%%%%%%%%%%%%%%%%%%%%%%%%%%%%%%%%%%%%%%%

\tableofcontents

%%%%%%%%%%%%%%%%%%%%%%%%%%%%%%%%%%%%%%%%%%%%%%%%

\section{Introduction}

\subsection{Nilmanifolds: definition and basic properties}

%%%%%%%%%%%%%%%%%%%%%%%%%%%%%%%%%%%%%%%%%%%%%%%%%%%%%%%%%%%%%%%%%%%%%%%%

A {\bf nilmanifold} is a compact manifold equipped with a transitive action
of a nilpotent Lie group. As shown by Malcev (\cite{_Malcev_}),
every nilmanifold can be obtained as a quotient of a nilpotent
Lie group $G$ by a discrete lattice $\Gamma$. Moreover, the
group $G$ can be obtained as so-called Malcev completion of $\Gamma$,
that is, as a product of exponents of formal logarithms of
elements $\Gamma$. Therefore, any nilmanifold is uniquely determined
by its fundamental group, which is a discrete nilpotent
torsion-free group, and any such group uniquely determines
a nilmanifold.

Since the work of Thurston  (\cite{_Thurston:Kodaira_}),
geometric structures on nilmanifolds are used
to provide many interesting examples (and counterexamples)
in complex and symplectic geometry. It was Thurston who realized
that the Kodaira surface (also known as a Kodaira-Thurston surface)
is symplectic, but does not admit any K\"ahler structure.
In this  way Thurston obtained a counterexample to
a result stated by H. Guggenheimer (\cite{_Guggenheimer:fail_})
in 1951. Guggenheimer claimed that the Hodge decomposition is true
for compact symplectic manifolds, but for symplectic
nilmanifolds this is usually false.

Before 1990-ies, a  ``complex nilmanifold'' meant a compact quotient
of a complex nilpotent Lie group by a discrete, co-compact subgroup.
The first non-trivial example is given by so-called Iwasawa manifold
(\cite{_Fernandes_Gray:Iwasawa_}) which is obtained as a quotient
of the 3-dimensonal Lie group of upper triangular
3 by 3 matrices by a discrete co-compact subgroup, for example the
group of upper triangular matrices with coefficients in $\Z[\1]$.

Starting from late 1980-ies,
a ``complex nilmanifold'' means a quotient
of a real nilpotent  Lie group equipped with a left-invariant
complex structure by the left action of a
discrete, co-compact subgroup (\cite{_CFG:_Frolicher_}).
This is the notion we are going to use in this paper.
This definition is much more general, indeed, left-invariant
complex structures are found on many even-dimensional
nilpotent Lie groups which are not complex.
The complex structure on a Kodaira surface
is one of such examples.

Complex structures on a nilmanifold have a very neat
algebraic characterization. Let $G$ be a real nilpotent
Lie group, and $\g$ is Lie algebra. By Newlander-Nirenberg
theorem, a complex structure on $G$ is the same as a
sub-bundle $T^{1,0}G \subset TG\otimes_\R \C$ such that
$[T^{1,0}G, T^{1,0}G]\subset T^{1,0}G$ and
$T^{1,0}G\oplus \overline{T^{1,0}G}=TG\otimes_\R \C$.
The left-invariant
sub-bundles in $T^{1,0}G$ are the same as subspaces
$W\subset \g\otimes_\R \C$, and the condition
$[T^{1,0}G, T^{1,0}G]\subset T^{1,0}G$ is equivalent
to $[W,W]\subset W$. Therefore, left-invariant
complex structures on $G$ are the same as
complex sub-algebras $\g^{1,0}\subset \g\otimes_\R \C$
satisfying $\g^{1,0}\oplus \overline{\g^{1,0}}=\g\otimes_\R \C$.

A real nilmanifold is obtained as an iterated fibration
with fibers which are compact tori.
It is natural to ask if any complex nilmanifold
can be obtained as an iterated fibration with
fibers which are complex tori. The answer is negative:
see e.g. \cite{_Rollenske:def_at_large_}.

However, a weaker statement is still true.
If we replace fibrations of nilmanifolds by
homomorphisms of their Lie algebras,
it is possible to construct a homomorphism
$\psi:\; \g \arrow \a$ to a complex abelian Lie algebra
compatible with a complex structure. Since
$\a$ is abelian, $\ker \psi$ necessarily contains
the commutator $[\g, \g]$. Since it is complex,
$\ker \psi$ contains $[\g, \g]+ I[\g, \g]$.

The quotient algebra $\g/[\g, \g]+ I[\g, \g]$
is called {\bf the algebra of holomorphic differentials
on $G$}, denoted by ${\goth H}^1(M)$. Its dimension is always positive
(\cite{_Sal_}).

In this paper, we study algebraic objects
(such as meromorphic functions) on complex
nilmanifolds. We prove the following theorem.

\hfill

%%%%%%%%%%%%%%%%%%%%%%%%%%%%%%%%%%%%%%%%%%%%%%%%%%%%%%
\theorem\label{_mero_constant_on_commu_Theorem_}
Let $M=G/\Gamma$ be a complex nilmanifold, and
$\Sigma$ be a  a foliation
obtained by left translates of $[\g, \g]+ I[\g, \g]$,
where $I$ is the complex structure operator, and
$\g=T_eG$ the Lie algebra of $G$.
Then all meromorphic functions on $M$
are constant on the leaves of
$\Sigma$.

\hfill

{\bf Proof:} See Subsection \ref{_ave_posi_Subsection_}. \endproof

%%%%%%%%%%%%%%%%%%%%%%%%%%%%%%%%%%%%%%%%%%%%%%%%%%%%%%%%%%%%%%%%%%%%%%%%
\subsection{Algebraic dimension and holomorphic differentials}
\label{_Kahler_rk_Subsection_}
%%%%%%%%%%%%%%%%%%%%%%%%%%%%%%%%%%%%%%%%%%%%%%%%%%%%%%%%%%%%%%%%%%%%%%%%

Recall that a positive closed (1,1)-current $T$ on a
complex manifold is said to have {\bf analytic
  singularities} (see \cite{bouck}) if locally
$T=\theta+dd^c\phi$ for a smooth form $\theta$ and
a plurisubharmonic function $\phi = c
\log(|f_1|^2+...+|f_n|^2)$ where $f_1,...f_n$ are analytic
functions and $c$  a constant.
Such currents have
decomposition into absolutely continuous and singular
part, where the absolutely continuous part is positive and
closed.

\smallskip

\definition  Let $M$ be a complex manifold.  The {\bf
  K\"ahler rank} $k(M)$ of $M$ is the maximal rank of the
absolutely continuous part of a positive, closed
(1,1)-current on $M$ with analytic singularities.

\hfill

\definition  The {\bf algebraic dimension} $a(M)$ of a complex manifold is
the transcendence degree of its field of meromorphic functions.

\hfill

Let $X$ be a complex manifold, and
$\phi:\; X \dashedrightarrow \C^n$ a meromorphic map
defined by generators of the field of
meromorphic functions. An algebraic reduction
of $X$ (\cite{_Campana:reduction_}, \cite{_Ueno_})
is a compactification of $\phi(X)$ in $\C P^n\supset
\C^n$. It is known to be a compact, algebraic
variety (\cite{_Campana:reduction_}, \cite{_Ueno_}).

We should note that the map $\phi$ is defined for more
general spaces $X$. For smooth manifolds we'll use the
following \cite[Definition-Theorem 6.5]{Peternel}.

\hfill

%%%%%%%%%%%%%%%%%%%%%%%%%%%%%%%%%%%%%%%%%%%%%%%%%%%%%%%%%%%%
\definition \label{_algebra_red_Definition_}
Let $M$ be a compact complex manifold. Then there
exists a smooth projective manifold $X$, a rational map
$\phi:\; M \dashedrightarrow X$ and a  diagram
\begin{diagram}[size=2em]
& & X' & &\\
& \ldTo^{a} & & \rdTo^{b} & \\
M & &  \rDashto_{\phi}  &&X \\
\end{diagram}
 where $X'$ is smooth and the top two arrows are proper
 holomorphic maps with $a$ a proper bimeromorphic
 modification, such that the 
corresponding fields of meromorphic functions 
coincide: $\Mer(M)=\Mer(X)$. We call
 the map $\phi:\; M \dashedrightarrow X$ {\bf algebraic
   reduction} of $M$.

\hfill

%%%%%%%%%%%%%%%%%%%%%%%%%%%%%%%%%%%%%%%%%%%%%%%%%%%%%%%%%%%%%%%
\definition\label{_induced_from_Kah_Definition_}
Let $\eta$ be a positive $(1,1)$-form on $X$.
The current $T_{\eta}$ is defined as $a_*b^* \eta$; since pushforward of a form
is a current, it is a current, and not a form.
Since $a$ is one-to-one everywhere, except on an analytic set $E\subset
X'$, the current $a_*b^* \eta$ is smooth outside of $E$.
Note also that the positivity and closedness are preserved,
as well as the rank in a general point.  We call $T_{\eta}$ the
{\bf current induced by $\eta$ on $M$}.
It is easy to check that $T_{\eta}$ has analytic
singularities if $\eta$ is closed and positive.

\hfill

%%%%%%%%%%%%%%%%%%%%%%%%%%%%%%%%%%%%%%%%%%%%%%%%
\claim  Let $M$ be a complex variety. Then
the algebraic dimension is bounded by  the K\"ahler rank:
\begin{equation}\label{_algebra_dime_Kahler_rank_Equation_}
a(M) \leq k(M).
\end{equation}

\hfill

{\bf Proof:} Let $\phi:\; M \dashedrightarrow X$ be the algebraic
reduction map. Pullback of a K\"ahler form from $X$ to $M$
is a current of rank $\dim X$ at all points where it is
absolutely continuous. \endproof

\hfill

We are going to estimate $a(M)$ in terms of holomorphic
differentials on $M$.

\hfill

%%%%%%%%%%%%%%%%%%%%%%%%%%%%%%%%%%%%%%%%%%%%%%%%
\definition  {\bf A holomorphic differential} on a compact complex manifold
is a closed, holomorphic 1-form.

\hfill

%%%%%%%%%%%%%%%%%%%%%%%%%%%%%%%%%%%%%%%%%%%%%%%%%%%%%%%%%%%%
\theorem\label{_main_alg_dim_Theorem_}
 Let $M$ be a complex nilmanifold, ${\goth H}^1(M)$
the space  of holomorphic
differentials on $M$, and $a(M)$ its
algebraic dimension. Then $$a(M) \leq \dim {\goth H}^1(M).$$

\hfill

{\bf Proof:} See Subsection \ref{_ave_posi_Subsection_}.
\endproof

\hfill

%%%%%%%%%%%%%%%%%%%%%%%%%%%%%%%%%%%%%%%%%%%%%%%%
\remark
The same estimate holds for complex parallelizable
manifolds; see \cite[Proposition 3.16.1]{_Winkelmann:paralleli_}.
Moreover, \ref{_mero_constant_on_commu_Theorem_} is also true
in this situation: all meromorphic functions are constant
on the fibers of the Albanese map.

\hfill

This result is implied by \ref{_algebra_red_Definition_}
and the following expression of the K\"ahler rank.

\hfill

%%%%%%%%%%%%%%%%%%%%%%%%%%%%%%%%%%%%%%%%%%%%%%%%
\theorem\label{_Kahler_rank_Theorem_}
Let $M$ be a complex nilmanifold, ${\goth H}^1(M)$ the space holomorphic
differentials on $M$, and $k(M)$ its  K\"ahler rank.
Then $k(M) =\dim {\goth H}^1(M)$.

\hfill

{\bf Proof:}
Consider the projection $\g \arrow \a$,
where $\a=\frac{\g}{[\g,\g]+I([\g,\g])}$.
Since $[\g,\g]+I([\g,\g])$ is $I$-invariant, $a$ has a complex structure and
this map is compatible with it. Consider the Chevalley differential
$d$ on the Lie algebras of $\g$ and $\a$.
Since $\a$ is an abelian algebra, any
2-form on $\a$ is closed (and gives a closed 2-form on the
corresponding Lie group). Taking
a positive definite Hermitian form, we obtain
a positive current of rank $\dim \a={\goth H}^1(M)$
on $M$. There are no currents with greater rank
by \ref{_current_folia_Proposition_}.
\endproof

\hfill

The same argument implies the following useful corollary.

\hfill

%%%%%%%%%%%%%%%%%%%%%%%%%%%%%%%%%%%%%%%%%%
\corollary\label{_factori_Kah_Corollary_}
Let $M$ be a complex nilmanifold, $\g$ the Lie algebra
of the corresponding Lie group, and ${\goth h}:= [\g,\g]+I([\g,\g])$
the algebra constructed as above. Denote by ${\goth h}_1$
a smallest $I$-invariant rational subspace of $\g$
containing ${\goth h}$. Let $T$ be a complex torus
obtained as quotient of $\g/{\goth h}_1$ by its integer lattice.
Consider the natural holomorphic projection $\Psi:\; M\arrow T$.
Then any meromorphic map to a K\"ahler
manifold is factorized through $\Psi$.

\hfill

{\bf Proof:}
Let $\psi:\; M \dashedrightarrow X$ be a meromorphic map
to a K\"ahler manifold $(X, \omega)$. For general $x\in X$, the
zero space of the positive closed current $\psi^*\omega$
contains ${\goth h}$, hence the fibers $F_x:=\psi^{-1}(x)$ are tangent
to ${\goth h}$. The smallest compact complex
subvariety of $M$ containing a leaf of the foliation
associated with ${\goth h}$ is the corresponding leaf of ${\goth h}_1$.
Passing to the closures of
the leaves of ${\goth h}$, we obtain that $F_x$ contain
leaves of ${\goth h}_1$. However, $T$ is the leaf space of ${\goth h}_1$.
\endproof

\hfill

\remark\label{rem}
For a general compact complex manifold $X$, {\bf Albanese variety} $\Alb(X)$ is defined as 
the quotient of the dual space of the space of
holomorphic differentials 
$H^0(X, d\mathcal{O})^*$ by the minimal closed complex
subgroup containing the image of $H^1(X,\mathbb{Z})$ under the map
$$ H^1(X,\mathbb{Z}) \rightarrow H^1(X,\mathbb{C}) \rightarrow H^0(X, d\mathcal{O})^*$$
(see \cite{_Rollenske:def_at_large_}). The Albanese map $\Alb: X\rightarrow \Alb(X)$ is give
n by integration along paths starting at a fixed point. It has the functorial property tha
t any map from $X$ onto tori is factored through the Albanese map (\cite{_Ueno_}).
In  \cite{_Rollenske:def_at_large_} (section 2.1) the Albanese variety for a complex nilma
nifold $M$ is described in terms of the space
$\mathfrak{h_1}$ as 
\[ \Alb(M) = \frac{ H^0(X, d\cal{O})^*/p(\mathfrak{h_1})}{\im(H^1(X,\mathbb{Z}) \rightarrow H^0(X, d\mathcal{O}))^*/p(\mathfrak{h_1})} = T,
\]
 where $\mathfrak{h_1}$ is the same  as in \ref{_factori_Kah_Corollary_}. 
Then we obtain that $T=\Alb(M)$ and $$a(M)=a(\Alb(M)).$$

%%%%%%%%%%%%%%%%%%%%%%%%%%%%%%%%%%%%%%%%%%%%%%%%

\section{The averaging formalism}
\label{_Averaging_Section_}

%%%%%%%%%%%%%%%%%%%%%%%%%%%%%%%%%%%%%%%%%%%%%%%%

Let  $M = \Gamma \backslash G$  be a compact nilmanifold
and $\nu$  a volume element on $M$ induced by a
the Haar measure on the Lie group $G$ \cite{Milnor}. After
a rescaling, we can suppose that  $M$ has volume
$1$. Notice that the Haar measure on $G$ is bi-invariant,
because $G$ admits a lattice, and any Lie group admitting
a lattice is unimodular.

Given any covariant $k$-tensor field $T : TM \times
\ldots \times T M  \rightarrow {\mathcal C}^{\infty}
(M)$  on the nilmanifold $M$,  one can define  a
covariant $k$-tensor $$T_{inv} : \frak g \times \ldots
\times \frak g \rightarrow \R$$  on the Lie algebra
$\frak g$  of $G$ by
$$
T_{inv} (x_1, \ldots, x_k) = \int_{p \in M} T_p (x_1 \vert_p, \ldots, x_k \vert_p) \nu,
$$
for every $x_1, \ldots, x_k \in \frak g$, where $x_l
\vert_p$ is the restriction of the
left-invariant vector field $X_l$ to $p$.
Clearly, $T_{inv} = T$  for any tensor field T coming
from a left-invariant one. In \cite{Belgun}   it is shown
that that  if $\alpha$ is a differential $k$-form on $M$,
then $(d\alpha )_{inv}  = d(\alpha_{inv})$.  Moreover,
$(\alpha_{inv}   \wedge \beta )_{inv}  = \alpha_{inv}
\wedge \beta_{inv}$, for every differential forms $\alpha$
and $\beta$ on $M$.

We call the map $\Av:\; (T^*)^{\otimes k}\arrow (\g^*)^{\otimes k}$,
$\Av(T):=T_{inv}$
{\bf averaging} on a nilmanifold.
The averaging defines a linear map $\tilde
\nu : \Omega^k (M)  \rightarrow \Lambda^k \frak g^*$,
given by  $\tilde \nu  (\alpha) = \alpha_{inv}$ for every
$k$-form $\alpha \in \Omega^k (M),$  which commutes with
the differentials.

Moreover, by Nomizu theorem
\cite{Nomizu}   $\tilde  \nu$ induces an isomorphism
$H^k(M) \rightarrow  H^k (\frak g)$ between the kth
cohomology groups for every  $k$. In particular, every
closed $k$-form $\alpha$ on M is cohomologous to the
invariant $k$-form $\alpha_{inv}$  obtained by the
averaging (see also \cite{_Ugarte_}) Indeed, by Nomizu Theorem $\alpha = \beta + d \gamma$, with $\beta$ invariant closed $k$-form. By using the averaging we have
$\alpha_{inv} = \beta + d \gamma_{inv}$ and so $\alpha$ is cohomologous to $\alpha_{inv}$.

%%%%%%%%%%%%%%%%%%%%%%%%%%%%%%%%%%%%%%%%%%%%%%%%

\section{Positive currents on nilmanifolds}

%%%%%%%%%%%%%%%%%%%%%%%%%%%%%%%%%%%%%%%%%%%%%%%%

%%%%%%%%%%%%%%%%%%%%%%%%%%%%%%%%%%%%%%%%%%%%%%%%
\subsection{Holomorphic differentials}
%%%%%%%%%%%%%%%%%%%%%%%%%%%%%%%%%%%%%%%%%%%%%%%%

Recall that {\bf holomorphic differentials} on a complex manifold are closed, holomorphic 1-forms.

\hfill

%%%%%%%%%%%%%%%%%%%%%%%%%%%%%%%%%%%%%%%%%%%%%%%%
\definition Let $M=\Gamma \backslash G$ be a nilmanifold.
A differential form on $M$ is called {\bf invariant}
if its pullback to $M$ is invariant with respect to the
left action of $G$ on itself.

\hfill

%%%%%%%%%%%%%%%%%%%%%%%%%%%%%%%%%%%%%%%%%%%%%%%%
\remark  Let $M=\Gamma \backslash G$ be a nilmanifold,  $\g$ the Lie algebra
of $G$. Clearly, invariant
differential forms are identified with $\Lambda^*(\g)$.
Moreover, they are preserved by de Rham differential,
which is identified with the Chevalley differential
on $\Lambda^*(\g)$.

\hfill

%%%%%%%%%%%%%%%%%%%%%%%%%%%%%%%%%%%%%%%%%%%%%%%%%%%%%%%%%%%%
\proposition  Let $(M,I)$ be a complex nilmanifold, and
$h$ a holomorphic differential. Then $h$ is
an invariant differential form.

\hfill

{\bf Proof:} Let  $\nu$  be  a volume element on $M$
induced by a bi-invariant one on the Lie group $G$  such
that $M$ has volume equal to $1$. A holomorphic
differential $h$ is cohomologous  to the invariant  form
$h_{inv}$  obtained by the averaging process. Since $I$ is
invariant, $h_{inv}$  has to be  of type $(1,0)$ and thus
$h =  h_{inv}.$
Indeed,   closed (1,0)-forms cannot be exact, because they are holomorphic,
hence (if exact) equal to differentials of a global holomorphic function.
 \endproof

\hfill

%%%%%%%%%%%%%%%%%%%%%%%%%%%%%%%%%%%%%%%%%%%%%%%%
\corollary\label{_dim_holo_diff_Corollary_}  Let $M=\Gamma \backslash G$ be a complex nilmanifold,
and $\g$ its Lie algebra, and ${\goth H}^1(M)$
the space of holomorphic differentials. Then
$${\goth H}^1(M)= \left(\frac{\g\otimes \C}{\g^1+ I(\g^1)}\right)^*,$$
where $\g^1=[\g,\g]$ denotes the commutator of $\g$.

\hfill

{\bf Proof:}  Let $h$ be a holomorphic differential. Since $h$ is invariant then it can be identified with an element of  $(\g \otimes \C)^*$. Moreover, $h = \alpha + i I \alpha$, with $\alpha \in \g^*$, $d \alpha=0$ and $d (I \alpha)=0$. By the conditions
$$
d \alpha (x, y) = - \alpha ([x, y])=0, \quad d (I \alpha) (x, y) = \alpha (I [x, y])=0,
$$
for every $x, y \in \g$,  we get $\alpha (\g^1) = \alpha (I \g^1) =0$.
 \endproof

%%%%%%%%%%%%%%%%%%%%%%%%%%%%%%%%%%%%%%%%%%%%%%%%
\subsection{Positive (1,1)-forms on a Lie algebra}
%%%%%%%%%%%%%%%%%%%%%%%%%%%%%%%%%%%%%%%%%%%%%%%%

Throughout this subsection, we fix a nilpotent
Lie algebra $\g$ with a complex structure
$I\in \End(\g)$ satisfying the integrability condition
$$
[\g^{1,0}, \g^{1,0}]\subset \g^{1,0}.$$

\hfill

\definition  A semipositive Hermitian
form on $(\g, I)$  is a real form $\eta\in \Lambda^2(\g^*)$
which is $I$-invariant (that is, of Hodge type
(1,1)) and satisfies $\eta(x, Ix)\geq 0$ for each
$x\in \g$. It is called {\bf positive definite Hermitian}
if this inequality is strict for all
$x\neq 0$.

\hfill

%%%%%%%%%%%%%%%%%%%%%%%%%%%%%%%%%%%%%%%%%%%%%%%%
\definition  A subalgebra $\a \subset \g$ is called
{\bf holomorphic} if $I(\a)=\a$ and
$[\g^{0,1}, \a^{1,0}]^{1,0}  \subset \a^{1,0}$.

\hfill

%%%%%%%%%%%%%%%%%%%%%%%%%%%%%%%%%%%%%%%%%%%%%%%%
\claim
Let $\a \subset \g$ be a vector subspace, and
$B:=\a \cdot G$ the corresponding left-invariant
sub-bundle in $TG$. Then
\begin{itemize}
\item $B$ is involutive (that is, Frobenius integrable)
iff $\a$ is a  Lie subalgebra of $\g$.
\item $B$ is a holomorphic sub-bundle iff $\a$ is a holomorphic subalgebra.
\end{itemize}

{\bf Proof:} 
  Let  $x, y \in \a$ and denote by the same letters
the corresponding  left-invariant vector fields. Clearly, $B$ is involutive if
and only if  $\a$ is a  Lie subalgebra of $\g$. Similarly we have
%$$
%[(x + i I  x)^R, (y - i I y)]^R =  - [x + i I x, y - i I y]^R
%$$
that $B$ is holomorphic if
$[x + i I x, y - i I y] \in {\frak a}^{1,0}$, for every $x \in \frak g$ and  $y \in \frak 
a.$ \endproof

\hfill

\remark
Note that $V = \g^{(1,0)} + \a^{(0,1)}$
is involutive iff $\a$ is holomorphic and
$V+\overline{V}=\g^c$. So $V$ is an \lq \lq elliptic
structure" in the terminology of \cite{Jac}, so by
  \cite{Jac} it defines a holomorphic
  foliation.

\endproof

\hfill

We also note the obvious

\hfill

\claim\label{intersection} If $V_1$ and $V_2$ are two
elliptic structures in terminology of \cite{Jac} on a
complex manifold, containing the $(1,0)$ tangent bundle,
then $V_1\cap V_2$ is also an elliptic structure.

\hfill

%%%%%%%%%%%%%%%%%%%%%%%%%%%%%%%%%%%%%%%%%%%%%%%%
\definition Let $\eta$ be a semipositive Hermitian
form on $(\g, I)$, and $N(\eta)$ the subspace in $\g$
consisting of all vectors $x$ such that $\eta(x, Ix)=0$.
Then $N(\eta)$ is called {\bf the null-space} of $\eta$.

\hfill

In general we have the following

\hfill

%%%%%%%%%%%%%%%%%%%%%%%%%%%%%%%%%%%%%%%%%%%%%%%%
\claim
The nullspace
$$
N  = \{ x \in \g \, \mid  \iota_x \eta =0 \}
$$
of a closed form  $\eta \in \Lambda^r \frak g^*$  is a  Lie subalgebra of $\frak g$.

\smallskip

{\bf Proof:} Take $x,y\in N$ and  arbitrary vectors $z_1, \ldots, z_{r -1} \in \frak g$.  Then, by Cartan's formula,  $d\eta (x,y,z_1, \ldots, z_{r -1})= \eta([x,y], z_1, \ldots z_{r-1})=0$, since
the rest of the terms vanish, because $x, y\in N$.
Therefore  $\eta([x,y], z_1, \ldots z_{r-1})=0$ for any $z_1, \ldots, z_{r -1} \in \frak g$, this means
that  $\iota_{[x,y]} \eta =0$, i.e $[x, y]\in N$.
\endproof

\hfill

%%%%%%%%%%%%%%%%%%%%%%%%%%%%%%%%%%%%%%%%%%%%%%%%
\theorem\label{_Kahler_quotient_Theorem_}
Let $\eta$ be a semipositive Hermitian
form on $(\g, I)$. Assume that its nullspace
$N(\eta)$ is a holomorphic subalgebra. Then
$N(\eta)$ contains $\g^1+ I\g^1$, where
$\g^1=[\g,\g]$.

\hfill

%
%
%\marginpar{WE SHOLD PUT A REAL PROOF HERE!}
%
%{\bf Proof. Step 1:} The leaf space of $G\cdot N(\eta)$
%
%is a complex variety equipped with a positive definite
%
%closed (1,1)-form.
%
%
%
%
%
%{\bf Step 2:} It is also a homogeneous space
%
%over a nilpotent Lie group, hence
%
%cannot be K\"ahler.

{\bf Proof:}
When the cohomology class of $\eta$ is rational,
as happens in most applications,
\ref{_Kahler_quotient_Theorem_} has a simple proof.
Since $[\eta]$ is rational, it can be represented
by a rational form $\eta_\Q\in \Lambda^2(\g)$.
Therefore, the leaves of $N(\eta)$
are rational Lie subalgebras in $\g$.
By Malcev's theorem, the leaves of  $N(\eta)$ are compact.
By construction, the leaf space $X$ of $N(\eta)$
is equipped with a transitive action by a
nilpotent Lie group, hence it is a nilmanifold.
Finally, $X$ inherits the complex structure
from $X$, and $\eta$ defines a K\"ahler metric
on $X$. However, a nilmanifold can be K\"ahler
only if its fundamental group is abelian (\cite{_BG_}).
Therefore, $N(\eta_\Q)$ contains $[\g, \g]$.

For general $\eta$, \ref{_Kahler_quotient_Theorem_}
has a different (more complicated) proof.

Since  $N(\eta) = \frak a $ is holomorphic, we have
\begin{equation} \label{holomcond}
[y + i I y, x - i I x]^{1,0} \in {\frak a}^{1,0},
\end{equation}
for every $x \in \frak a$ and for every $y \in \frak g$. By a direct computation we obtain
$$
\begin{array}{lcl}
[y + i I y, x - i I x]^{1,0} &=& ([y, x] + [Iy, Ix] + I [Iy, x] - I [y, Ix])\\[3pt]
&& - i I ([y, x] + [Iy, Ix] + I [Iy, x] - I [y, Ix]).
\end{array}
$$
Therefore,  by the condition \eqref{holomcond} we get
\begin{equation}\label{holom2}
[y, x] + [Iy, Ix] + I [Iy, x] - I [y, Ix] \in \frak a , \quad \forall x \in \frak a,  \, \forall y \in \frak g.
\end{equation}
By  using the integrability condition
$$
[Iy, Ix] = [y,x] + I[Iy, x] + I [y, Ix]
$$
we have
$$
I[y, Ix] = [Iy, Ix] - [y, x] - I [Iy, x]
$$
and therefore the condition \eqref{holom2} becomes
$$
2 ([y, x] + I [Iy, x])  \in \frak a , \quad \forall x \in \frak a,  \, \forall y \in \frak g,
$$
i.e.
$$
\eta ([y, x], z) = - \eta (I [Iy, x], z), \quad  \forall x \in \frak a,  \, \forall y, z \in \frak g.
$$
Therefore
$$
\eta ([y, x], I [y, x]) = - \eta (  I [Iy, x],  I [y, x]) =  - \eta ([x, Iy], [x, y]).
$$
By $d\eta =0$, one gets
$$
 \eta ([x, y],w) = \eta ([x, w], y),
 $$ for every $x \in \frak a, y, z \in \frak g.$
Thus
$$
\eta ([x, Iy], [x, y]) = \eta (ad_x^2 (y), Iy), \quad \forall x \in \frak a, \, \forall y \in \frak g
$$
and consequently
\begin{equation} \label{condad2}
\eta ([y, x], I [y, x])  = - \eta (ad_x^2 (y), Iy), \quad \forall x \in \frak a, \, \forall y \in \frak g.
\end{equation}
 By using \eqref{condad2}, it is possible to show that
 $\frak a$ is an ideal of $\frak g$,  i.e. that  $[y, x]
 \in \frak a$, for every $x \in \frak a$ and $y \in \frak g$.

Since $\eta$ is a semipositive (1,1)-form and $\goth a$ is
its null-space, the relation $\eta([y, x], I [y, x])=0$ implies that
$[x,y]\in \goth a$. Therefore, by \eqref{condad2}, in
order to prove that $[y, x] \in \frak a$,  for every $x
\in \frak a$ and for every $y \in \frak g$, it is
sufficient to show that $[x, [x, y]]\in \goth a$
for any $x\in \goth a$. This would follow if we prove that
\begin{equation}\label{_a_with_commu_Equation_}
[\frak a, {\frak g}^1] \subset  \frak a.
\end{equation}
Since $\frak g$ is nilpotent there exists $s$ such that $\frak g^s = \{ 0 \}$ and $\frak g^{s-1} \neq  \{ 0 \}$ and we have the descending series of ideals
$$
\frak g = {\frak g}^0 \supset   {\frak g}^1   \supset  \ldots \supset  {\frak g}^i \supset {\frak g}^{i + 1} \supset   \ldots \supset  {\frak g}^{s - 1} \supset {\frak g}^s = \{ 0\}.
$$
Now  we can prove that $[\frak a, {\frak g}^1] \subset
\frak a$   by induction on  $i$   in the following way:
by using \eqref{condad2} we can show that

\smallskip

(A) if the condition $[\frak a, {\frak g}^{i +1}] \subset \frak a$   holds, then  the condition  $[\frak a, {\frak g}^i] \subset \frak a$ holds.

\smallskip

Since $\frak g$ is nilpotent there exists
$s$ such that $\frak g^s = \{ 0 \}$ and $\frak g^{s-1}  \neq \{ 0 \}$.

At  the first step  $i = {s-1}$ we have that   (A) holds.
So by induction  we  obtain that (A)  holds  for $i = 0$.
Consequently,  $\frak a = N(\eta)$ is an ideal of $\frak
g$ and  $\eta$ induces a K\"ahler form on  the nilpotent
Lie algebra ${\frak g} / {\frak a}$.  By \cite{_BG_}, the K\"ahler
nilpotent Lie algebra ${\frak g} / {\frak a}$ has to be
abelian. Therefore ${\frak g}^1 \subset  \frak a$.
Since $\eta$ is (1,1)-form, its null-space
$\frak a$ is $I$-invariant,  hence $\frak a$ contains
$\frak g^1  + I \frak g^1$.
\endproof

%%%%%%%%%%%%%%%%%%%%%%%%%%%%%%%%%%%%%%%%%%%%%%%%

\subsection{Averaging the positive currents}
\label{_ave_posi_Subsection_}

%%%%%%%%%%%%%%%%%%%%%%%%%%%%%%%%%%%%%%%%%%%%%%%%

The following result directly follows from the averaging.

\hfill

%%%%%%%%%%%%%%%%%%%%%%%%%%%%%%%%%%%%%%%%%%%%%%%%
\proposition \label{_ave_posi_Proposition_}
Let $M=\Gamma \backslash G$ be a compact quotient of a
unimodular Lie group $G$ by a lattice $\Gamma$ and $I$
the complex structure on $M$ obtained from
an invariant complex structure on $G$. Let $T_{\eta}$ be the
positive, closed (1,1)-current induced by the algebraic
reduction $\phi:M\rightarrow X$ from some K\"ahler form
$\eta$ on $X$ (\ref{_induced_from_Kah_Definition_}). If $\Av(T)$ is its average, then $\Av(T)$
is a semipositive, closed, $G$-invariant differential form,
and its rank is no less than the rank of  the absolutely
continuous part of $T_{\eta}$.

\hfill

{\bf Proof:} If $X$ and $Y$ are left-invariant vector fields on
$M$, then $T_{\eta}(X,Y)$ is a measurable function
when we consider $T_{\eta}$ as a form with distributional
coefficients in local coordinates. So $\Av(T)$ is well
defined as in Section 2. Then $\Av(T)$ is a closed
invariant form of type (1,1) and the only thing to check
is the statement about its rank. By the definition it
follows that $\Av(T)(X,IX)=0 \Leftrightarrow
T_p(X|_p,IX|_p)=0$ for almost all $p\in M$. So $X$ is in
the kernel of $\Av(T)$ only if it is in the kernel of $T_p$
for almost all $p$. \endproof

\hfill

\remark
As a corollary we obtain that if such space
admits a K\"ahler current, it is K\"ahler. In particular
from \cite{Demailly-Paun} it follows that such spaces are
never in Fujiki's class $\mathcal{C}$. Note that the proof
of this fact in \cite{Demailly-Paun} uses also the
K\"ahler current arising from the pull-back of a K\"ahler
form.

\hfill

%%%%%%%%%%%%%%%%%%%%%%%%%%%%%%%%%%%%%%%%%%%%%%%%%%%%%%%%%%%%%
\proposition\label{_current_folia_Proposition_}
Let $T$ be a positive, closed (1,1)-current on a
nilmanifold $M=G/\Gamma$, and ${\cal F}$ the null-space foliation of
its absolutely continuous part. Then the sub-bunlde
associated with ${\cal F}$ contains a homogeneous sub-bundle
$\Sigma$ obtained by left translates of $\g^1+ I\g^1$,
where $\g^1=[\g,\g]$, and $\g$ is the Lie algebra of $G$.

\hfill

{\bf Proof:}
Let $\Av$ be the averaging map defined in
Section \ref{_Averaging_Section_}.  The nullspace of the form
$\Av(T)$ is contained in the intersection of all left translates
of ${\cal F}$, hence by  \ref{intersection}
it is also holomorphic.  Then
\ref{_Kahler_quotient_Theorem_}
implies that $N(\Av(T))$
contains $\g^1+ I\g^1$.
\endproof

\hfill

{\bf Proof of \ref{_mero_constant_on_commu_Theorem_}:}
Let $M \arrow X$ be the algebraic reduction
map (\ref{_algebra_red_Definition_}), and $\eta$ the
pullback of the K\"ahler form on $X$.
Averaging (\ref{_ave_posi_Proposition_})
transforms $\eta$ into an invariant, closed, semipositive form.
Then $\eta$ vanishes on $\Sigma$ by \ref{_current_folia_Proposition_}.
\endproof

\hfill

{\bf Proof of \ref{_main_alg_dim_Theorem_}:}
Now we can prove \ref{_main_alg_dim_Theorem_}.
Let $M = \Gamma \backslash G$ be a nilmanifold, and $\phi:\; M \dashedrightarrow X$
the algebraic reduction map. The pullback $\phi^* \omega_X$  of a K\"ahler
form $\omega_X$ is a current $T$ on $M$ (\ref{_induced_from_Kah_Definition_}). By
\ref{_current_folia_Proposition_}, the rank of its absolutely continuous part is no
greater that
$$\dim \frac{\g}{\g^1+ I\g^1}=\dim {\goth H}^1(M)$$
(the latter equality follows from
\ref{_dim_holo_diff_Corollary_}).
\endproof

\section{Examples}

All 2-dimensional compact complex nilmanifolds are
classified and correspond to tori and primary Kodaira surfaces. 
Their algebraic dimension is known. In this
section we'll consider the algebraic dimension of the
complex nilmanifolds in dimension 3 and note that for other
complex homogeneous spaces the inequality $a(M)\leq \frak
H^1 (M)$ may not hold.

Many nilmanifolds admit holomorphic fibrations and we'll need the following:

\smallskip

\remark
\label{remholfib} 
In general (see \cite [Theorem 3.8]{_Ueno_}), if  a
complex manifold $M$ is the total  space of a holomorphic
fibration $\pi: M \rightarrow B$ we always have the
inequality
$$
a(M) \geq a(B).
$$

\smallskip

\subsection{Algebraic dimension of complex 2-tori}

Following \cite{_Birkenhake_Lange_} we have the following description of the algebraic dimension of the complex 2-tori.

Let $T^4$ be the tori defined as ${\mathbb R}^4/{\mathbb Z}^4$ where ${\mathbb Z}^4$ is the standard lattice in ${\mathbb R}^4$. Let
$J\in \End({\mathbb R}^4),  J = \left (  \begin{array}{lcl}
A & B\\
C& D
\end{array}
\right )
$
be a complex structure with   $A, B, C, D$ $2\times
2$-blocks and $B$ nondegenerate. From
\cite{_Birkenhake_Lange_} (p.10) we can identify what is
the period lattice of the complex tori with structure
$J$. If $X= {\mathbb C}^2/(\tau, Id_2) {\mathbb Z}^4$ is a
complex tori defined by a complex $2\times 2$ matrix
$\tau$, then the complex structure $J_{\tau}$ on $T^4$
such that $X \cong (T^4, J_{\tau})$ as complex manifold is
given by  
$$
J_{\tau} = \left (  \begin{array}{lcl}
y^{-1}x & y^{-1}\\
-y-xy^{-1}x & -xy^{-1}
\end{array}
\right ),
$$
 where $x$ and $y$ are the real and imaginary parts of
 $\tau$. Reversing the construction gives that for $J$ as
 above, $$\tau_J= B^{-1}A+iB^{-1}.$$ We also need the
 relation to complex structures, defined in terms of a
 basis of $(1,0)$-forms. If $J_0$ is a fixed complex
 structure and $\omega_j = e_j + \1 J_0 e_j, j=1,2$ is a basis
 of (linear) $(1,0)$-forms for $J_0$, we define another
 complex structure $J$ as
\begin{equation}\label{formcxstr}
\begin{array}{lll}
\alpha_1 &=& \omega_1 + a\overline{\omega_1}+b\overline{\omega_2}\\
\alpha_2 &=& \omega_2 + c\overline{\omega_1}+d\overline{\omega_2}
\end{array}
\end{equation}
being the basis of $(1,0)$-forms of $J$. If $X = \left (  \begin{array}{lcl}
a & b\\
c & d
\end{array}
\right ) = X_1 + iX_2$, then the relation between $X$ and the matrix representing $J$ in the basis $(e_i, J_0e_i)$ is given by $$J = \left (  \begin{array}{lcl}
Id+X_1 & X_2\\
X_2 & Id-X_1
\end{array}
\right )^{-1}\left (  \begin{array}{lcl}
0 & Id\\
-Id & 0
\end{array}
\right )\left(\begin{array}{lcl}
Id+X_1 & X_2\\
X_2 & Id-X_1
\end{array}\right)$$
We'll use the explicit form of $J$ and $\tau_J$ when $X_1 = \left (  \begin{array}{lcl}
0 & a\\
0 & 0
\end{array}
\right )
$ and $X_2 = \left (  \begin{array}{lcl}
0 & b\\
0 & 0
\end{array}
\right )$. Direct calculation (using the fact that $X_1^2=X_2^2=X_1X_2=X_2X_1=0)$ gives $J = \left (  \begin{array}{lcl}
2X_2 & Id-2X_1\\
-Id-2X_1 & -2X_2
\end{array}
\right )$ and $$\tau_J  = \left (  \begin{array}{lcl}
i & 2a+ 2bi\\
0 & i
\end{array}
\right ) = i ID + 2X.$$

To determine the algebraic dimension of $(T^4, J)$, we need first the Neron-Severi group $NS(J)$ of $J$. Let $\tau_{ij}$ are the components of $\tau_J$ and $$E = \left (  \begin{array}{llll}
0&a&b&c\\
-a&0&d&e\\
-b&-d&0&f\\
-c&-e&-f&0
\end{array}
\right )\in M_4({\mathbb Z})$$ be an integral matrix. Then $NS(J) =\{ E\in M_4(Z)| a+d\tau_{11}-b\tau_{12} +f\tau_{21}-c\tau_{22} + e  det(\tau) =0\}$. With these notations in mind, the algebraic dimension of $(T^4,J)$ is determined by

$$
a(J) = \frac{1}{2}\max\{ \rank (J^TE) | E\in NS(J), J^T E \geq 0\},$$
where the superscript $T$ means transposition. Note that not all complex structures are described in this way - we have the non-degeneracy condition on $B$ which is required for $(\tau, Id_2)$ to be a period matrix. It is well known that $a(J)$ could be any of 0,1 or 2.

The integrality condition leads to the fact that
generically $a(J)=0$. For $a(J)\geq 0$ from
\cite{_Birkenhake_Lange_} p.59, we know that
$a(J)=1$ exactly when the torus admits a period matrix
$(\tau, Id_2)$ with 
\[ \tau = \left (  \begin{array}{lcl}
\tau_1 & \alpha\\
0& \tau_2
\end{array}
\right )
\] with 
\[ \alpha\notin (\tau_1,1)M_2(\mathbb{Q})\left (  \begin{array}{l}
1\\
\tau_2
\end{array}
\right ),
\] where $M_2(\mathbb{Q})$ is the set of $2\times 2$-matrices.
In particular in (\ref{formcxstr}), when $X=\left (  \begin{array}{lcl}
0 & \sqrt{2}-i\sqrt{3}\\
0 & 0
\end{array}
\right )$, the algebraic dimension of $(T^4, J)$ is one.

\subsection{3-dimensional complex nilmanifolds}

\begin{definition}   Let  $\frak g$  be a nilpotent Lie algebra. A rational structure for  $\frak g$  is a
subalgebra  $\frak g_{\mathbb Q}$  defined over  $\mathbb Q$  such that  $\frak g \cong  \frak g_{\mathbb Q} \otimes \R$.  A subalgebra  $\frak h$ of $\frak g$  is said to be rational with respect to a given rational structure $\frak g_{\mathbb Q}$  if  $\frak h_{\mathbb Q} :=  \frak h  \frak g_{\mathbb Q} $ is a rational structure for  $\frak h$.
\end{definition}

It follows from a result of Malcev \cite{_Malcev_} that
$\Gamma \backslash G$ is compact,  where G is a simply
connected k-step nilpotent Lie group admitting a basis of
left invariant 1-forms for which the coefficients in the
structure equations are rational numbers, and $\Gamma$ is
a lattice in $G$ of maximal rank (i.e., a discrete uniform
subgroup, cf. \cite{_Raghunatan_}). Such a lattice
$\Gamma$ exists in $G$ if and only if the Lie algebra
$\frak g$  of  $G$  has a rational structure.   Indeed, If
$\Gamma$  is a lattice  of  $G$ then its associated
rational structure is given by the $\mathbb Q$-span of
$\log \Gamma$.

An invariant  complex structure $J$  on  a nilmanifold $\Gamma \backslash G$  is called rational if it is compatible with the rational structure of $G$ , i.e. $J(\frak g_{\mathbb Q})  \subseteq   \frak g_{\mathbb Q}.$

\smallskip

\remark
Let $M = \Gamma \backslash G$ be a complex nilmanifold of complex dimension $n$  endowed with an  invariant  rational complex structure $J$.  Consider the surjective homomorphism $\frak g \rightarrow \frak g/\frak g^1_J$, where  $\frak g^1_J= \frak g^1 + J \frak g^1$.
Let  $G$, $G^1_J$  and $K$   be the simply connected Lie groups  respectively with Lie algebra $\frak g$, $\frak g^1_J$  and  $\frak g/\frak g^1_J$, then we have the surjective homomorphism
$$
p: G \rightarrow K,
$$
with $K$ abelian.
Since $J$ is rational then, by \cite{_Console_Fino_}  $\frak g^1_J$ is a rational subalgebra of $\frak g$.
Then $\Gamma^1 := \Gamma \cap  G^1_J$  is a uniform
discrete subgroup of $G^1_J$ \cite{_Corwin_Greenlaf_},
Theorem 5.1.11. By \cite[Lemma 5.1.4
  (a)]{_Corwin_Greenlaf_}, $p(\Gamma)$
is a uniform discrete subgroup of $K$ (i.e. $p(\Gamma) \backslash  K$) is compact, cf. \cite{_Raghunatan_}). By Lemma 2 in \cite{_Console_Fino_} the map
$$
\tilde p: \Gamma \backslash G \rightarrow   p(\Gamma) \backslash  K
$$
is a holomorphic fibre bundle. Moreover, since $K$ is
abelian, $p(\Gamma) \backslash  K$ is a complex torus
$\mathbb T$ of complex dimension $\frak H^1 (M) = n -
\dim_{\C} \frak g^1_J$. Therefore,  if $\frak  H^1 (M)
=1$, the torus $\mathbb T$ is algebraic, and by \ref{remholfib}
and  \ref{_main_alg_dim_Theorem_} we have
$
a (M) = 1.
$

\hfill

\begin{remark} Let $M = \Gamma \backslash G$ be a complex nilmanifold of complex dimension $n$  endowed with  an invariant complex structure $J$ such that $J \frak g^1 =  \frak g^1$, then $M$ is the total  space of a holomorphic fibration $\pi: M \rightarrow  \mathbb T$, with  $\mathbb T$  a complex torus of complex dimension $\frak H^1 (M) =   n - \dim_{\C}  \frak g^1$.
Therefore, if $\mathbb T$  is algebraic, i.e. $a(\mathbb
T) =  \frak H^1 (M)$, we have 
by \ref{remholfib} and \ref{_main_alg_dim_Theorem_}  
we have  $\frak a (B)= \frak a (M)  =  \frak H^1 (M).$

Note that  if  $J$ is bi-invariant, i.e.  if $M = \Gamma
\backslash G$ is  complex  parallelalizable then   $J
\frak g^1 =  \frak g^1$.    For  a general result on the
algebraic dimension of complex parallelalizable
solvmanifolds see Theorem 2 and its Corollary in
\cite{_Sakane_}.
\end{remark}

\hfill

We will apply the previous remarks to  complex  nilmanifolds  of complex dimension $3$.

\hfill

\example
Let $J$ be a complex structure on a  real $6$-dimensional
nilpotent Lie algebra. For the notion of ``nilpotent complex
structure'' on a nilmanifold, please see \cite{_Ugarte_}.
By \cite{_Ugarte_} the complex structure  $J$ is either nilpotent or non-nilpotent and
\begin{enumerate}
 \item[(a)]  If $J$  is non nilpotent, then there is a basis  of $(1,0)$-forms $(\omega^1, \omega^2, \omega^3)$ such that
$$
\left \{  \begin{array}{l}
 d \omega^1 =0,\\[3pt]
 d \omega^2 = E\,  \omega^1 \wedge \omega^3 + \omega^1  \wedge \overline \omega^3,\\[3pt]
 d \omega^3 = A \,  \omega^1 \wedge \overline \omega^1 + i b \, \omega^1 \wedge \overline \omega^2 - i b  \overline E \, \omega^2 \wedge \overline \omega^1,
 \end{array}
 \right.
 $$
 where $A, E \in \C$ with $| E | = 1$ and $b \in \R - \{ 0 \}$.
\item[(b)] If $J$ is nilpotent, then there is a basis of $(1,0)$-forms $(\omega^1, \omega^2, \omega^3)$  satisfying
$$
\left \{  \begin{array}{lcl}
 d \omega^1 &=& 0,\\[3pt]
d \omega^2 &= &\epsilon \omega^1 \wedge \overline \omega^1,\\[3pt]
d \omega^3 &= & \rho \, \omega^1 \wedge \omega^2 + (1 - \epsilon) A \, \omega^1 \wedge \overline \omega^1 + B \,  \omega^1 \wedge \overline \omega^2\\[2pt]
&&  + C \,  \omega^2 \wedge \overline \omega^1 + (1 - \epsilon) D \,  \omega^2 \wedge \overline \omega^2,
\end{array}
\right.
$$
where $A, B, C, D \in \C$ and $\epsilon, \rho \in \{ 0, 1 \}$.
 \end{enumerate}

 Suppose that the real and imaginary parts of the complex structure  equations constants are rational, then $G$ admits a lattice $\Gamma$. Let $M = \Gamma \backslash G$ be the compact quotient endowed with the induced invariant complex structure $J$.

  In the case $(a)$ we have that $\frak H^1 (M) =1$. In the case $(b)$ we have the following cases:
  \begin{enumerate}
\item[(b1)] $ \frak H^1 (M) =1$ if $\epsilon = 1$ and $\rho^2 + |B|^2  +   |C|^2 \neq 0$

\item[(b2)] $ \frak H^1 (M) =2$ if $\epsilon = 0$

\item[(b3)] $ \frak H^1 (M) =2$  if if $\epsilon = 1$ and $\rho = B=    C  =0$

 \end{enumerate}

 Therefore,  in the cases $(a)$ and $(b1)$,  since   $J$
 is rational,   by \ref{_main_alg_dim_Theorem_} 
and previous remarks, we have $a(M) =  \frak H^1 (M)
=1$. In the case $(b3)$, $G$ is the direct  product  of
the real $3$-dimensional Heisenberg  group by $\mathbb
R^3$.

  In the case $(b2)$, $ \rho^2 + |B|^2 + |C|^2 + |D^2| = 0$, then $G$ is the direct  product  of the real $3$-dimensional Heisenberg  group by $\mathbb R^3$.
  If  $ \rho^2 + |B|^2 + |C|^2 + |D^2| \neq 0$, then  $J \frak g^1 = \frak g^1$  is a rational subalgebra  of complex dimension $1$ and $M$ is the total space of a holomorphic fibre bundle over a complex torus $\mathbb T$  of complex dimension $2$.   Therefore,   if $\mathbb T$ is algebraic then by previous remarks we have  $a( M) =  \frak H^1 (M) = 2$.

  An explicit example of the case (b2) is given by the well known Iwasawa manifold $M$. The Iwasawa manifold $M$  is defined as the quotient $\Gamma \backslash G$, where
  $$
  G = \left \{ \left (  \begin{array}{ccc} 1&z_1&z_3\\ 0&1&z_2\\ 0&0&1 \end{array} \right )    \,  \mid  \, z_i \in \mathbb C \right  \}
   $$
   is the complex Heisenberg group and  $\Gamma$  is the lattice defined by taking $z_i$ to be Gaussian
integers, acting by left multiplication.  The 1-forms
$$
   \omega^1 = dz_1, \omega^2 = dz_2,  \omega^3  =  - dz_3  + z_1 d z_2
   $$
   are left-invariant on  $G$.  Define a  r basis $(e^1, . . . , e^6)$ of real 1-forms by setting
   $$
\omega^1 = e^i + i e^2,  \omega^2 = e^3 + i e^4,  \omega^6 = e^5 + i e^6.
   $$
   These 1-forms are pullbacks of corresponding 1-forms on the  compact quotient $M$, which we denote
by the same symbols and they  satisfy the structure equations
$$
   \begin{array}{l}
   d e^j =0, \,  j = 1,2,3,4, \\
   d e^5 = e^1 \wedge e^3 - e^2 \wedge e^4,\\
   d e^6 =  e^1 \wedge e^4 + e^2 \wedge e^3.
   \end{array}
   $$
 The Iwasawa manifold $M$ is the total space of a
 principal $T^2$-bundle over the real torus $T^4$. The
 mapping  $p:M  \rightarrow  T^4$ is induced from  the
 projection $(z_1, z_2, z_3) \mapsto (z_1, z_2)$ and
the space of invariant 1-forms annihilating the fibres of
p is  given by $\Span\langle e^1, e^2, e^3, e^4\rangle$

%    By \cite[Theorem 1.1]{_Ketsetzis_Salamon_}  let $J$ be
%    any invariant complex structure $J$ on $M$.

Then $p:M
    \rightarrow  T^4$ induces a complex
structure  $\hat J$ on  the real 4-dimensional torus
$T^4$  such that  $p: (M, J)  \rightarrow  (T^4,  \hat J)$
is holomorphic.

\hfill

\claim\label{iwasawaalgdim} For the invariant
  complex structures on the Iwasawa manifold,
$M$, $a(M) = a(T^4,\hat J)$

\hfill

{\bf Proof:} From \ref{_mero_constant_on_commu_Theorem_},
any meromorphic function is constant on the fibers of the
projection $M \arrow (T^4,\hat J)$. This implies that
$a(M) = a(T^4,\hat J)$. \endproof

\hfill

Now consider the possible algebraic dimension of $M$. First note that
not all complex structures on the base $T^4$ arise as $\hat J$. A description of  the set of such $J$ is given in [KS] and it is known that it has 4 components. We'll use one of them to see that that invariant structures on the Iwasawa nilmanifold can have algebraic dimension 0, 1 and 2. If we start with the canonical structure on $T^4$ which corresponds to the standard lattice in $R^4$ it gives the case of algebraic dimension 2. From \cite{_Ketsetzis_Salamon_} formula (10) we know that the matrix $X = \left(\begin{array}{ll}
  0& \sqrt{2}-i\sqrt{3}\\
  0& 0\\
  \end{array}\right)$ corresponds to a complex structure on $T^4$  which is $\hat J$ for some $J$ on $M$. On the other side, the period matrix $\tau$, as explained above is $iId+2X$ and has $\tau = x+iy$ with $x=
 \left(\begin{array}{ll}
  0& 2\sqrt{2}\\
  0& 0\\
  \end{array}\right)$ and $y=\left(\begin{array}{ll}
  1& -2\sqrt{3}\\
  0& 1\\
  \end{array}\right)$. So for such $\tau$ the algebraic
 dimension of the base is 1. In particular we have strict
 inequalities in \ref{_main_alg_dim_Theorem_}. Again, the base generically
 has algebraic dimension 0, which leads to vanishing of
 the algebraic dimension of $J$ on the Iwasawa
 nilmanifold.

%It remains only to check the \ref{iwasawaalgdim} above.
%It follows from the fact that every complex subvariety of
%codimension 1  on $M$ is a pull-back of a subvariety on
%$(T^4, \hat J)$ under $p$. To see this, assume that there
%is a subvariety $V\subset M$ which is transversal to a
%fiber of $p$ at at least one point. Than the form
%$\omega_3\wedge \overline{\omega_3}$ is non-vanishing and
%positive at at least one point f $V$. In particular
%$\int_V\omega_3\wedge\overline{\omega_3} >0$. But is
%exact on all $M$, so its integral over any oriented
%4-dimensional subspace vanishes - a contradiction. Then
%$V$ is fibered over $T^4$ and we have the same
%description of every divisor on $M$. Then it follows also
%$a(M)\leq a(T^4,\hat J)$, so combined with Remark
%\ref{remholfib} the Claim follows.

\smallskip

{\bf Example} (Compact Lie groups)  It is well known that
every even-dimensio\- nal compact Lie group $G$ admits an
invariant complex structure (\cite{Samelson}). The
construction uses the structure theory for semisimple Lie
algebras and provides a holomorphic fibration
$G\rightarrow Fl$ to the complete flag manifold $Fl =
G/T$, where $T$ is a maximal torus in $G$. The manifold
$Fl$ is algebraic and its algebraic dimension is equal to
its complex dimension. On the other side, $\frak g^{ss} =
\frak g'$, so the space of holomorphic differentials is
trivial. Hence by  \ref{remholfib}, $a(M)\geq dim(Fl)
>\frak H^1(M)=0$. Similarly, for non-K\"ahler compact
complex homogeneous spaces $G/H$ with $G$ compact, the
inequality $a(M)\leq \frak H^1 (M)$ does not hold in
general.

\hfill

{\bf Acknowledgements:} The work on this project started
when the first two authors visited the HSE in Moscow. They
are grateful for the hospitality and stimulating
environment at the "Bogomolov Lab" there. Part of the work
was done while the second named author visited University di Torino
(Torino), Max Plank Institute for Mathematics (Bonn) and
Institute of Mathematics at the Bulgarian Academy of
Sciences. He thanks all institutions for the
hospitality. The visits were partially funded by grants
from G.N.S.A.G.A. and MPIM. The third named author is 
thankful to Alexandra Victorova and Alexandra Skripchenko for
their ideas and many stimulating discussions of the subject.
We also thank S. Rollenske for his useful comments
on the paper.

{\small

\noindent
{\sc Anna Fino}\\
{\sc Dipartimento di Matematica G. Peano\\
Universit\'a di Torino\\
via  Carlo Alberto 10, 10123 Torino, Italy\\
\tt annamaria.fino@unito.it}\\

\noindent
{\sc Gueo Grantcharov\\
{\sc Department of Mathematics and Statistics\\
Florida International University\\
Miami Florida, 33199, USA}\\
\tt grantchg@fiu.edu}\\

\noindent {\sc Misha Verbitsky\\
{\sc Laboratory of Algebraic Geometry,\\
National Research University Higher School of Economics,\\
Department of Mathematics, 7 Vavilova Str. Moscow, Russia,}\\
\tt  verbit@mccme.ru}, also: \\
{\sc Universit\'e Libre de Bruxelles, CP 218,\\
Bd du Triomphe, 1050 Brussels, Belgium}

}


{\small
\begin{thebibliography}{AV1}



\bibitem[BD]{_BD:abelian_}
M.L. Barberis, I. Dotti,  {\em Abelian complex structures on
solvable Lie algebras}, J. Lie Theory {\bf 14}(1) (2004), 25--34.


\bibitem[BDV]{_BDV:nilmanifolds_}
M.  L. Barberis, I.  G. Dotti, M. Verbitsky, {\em Canonical bundles of complex nilmanifolds,
with applications to hypercomplex geometry}, Math. Res. Lett. {\bf 16} (2009), 331--347.

\bibitem[B]{Belgun}  F. Belgun, {\em On the metric structure of non-K\"ahler complex surfaces}, Math. Ann. {\bf 317} (2000), 1--40.

\bibitem[BG]{_BG_}
C. Benson, C.S. Gordon, {\em K\"ahler and symplectic structures on
nilmanifolds}, Topology {\bf 27}(4) (1988) 513--518.

\bibitem[BL]{_Birkenhake_Lange_} C. Birkenhacke, H. Lange,   Complex Tori, Progress in Mathematics 177, Birkhauser, 1999.

\bibitem[Bouck]{bouck} S. Boucksom. {\em On the volume of a line bundle} Intern. J. Math. {\bf 13}, no 10 (2002) 1043 -- 1063.

\bibitem[Ca81] {_Campana:reduction_}
F. Campana,
{\em Cor\'eduction alg\'ebrique d'un espace analytique faiblement
k\"ahl\'erien compact}, Invent. Math. 63 (1981) 187-223.


\bibitem[CF]{_Console_Fino_}
S. Console, A. Fino, {\em Dolbeault cohomology of compact
nilmanifolds}, Transformation Groups {\bf 6} (2) (2001), 111--124.

\bibitem[CFP]{_Console_Fino_Poon_}
S. Console, A. Fino, Y.S. Poon,  {\em Stability of abelian complex
structures,} Internat. J. Math. {\bf 17} (2006), no. 4, 401--416.



\bibitem[CFG1]{_CFG:_symple_}
L.A. Cordero, M. Fern\'andez, A. Gray,  {\em Symplectic manifolds
with no K\"ahler structure}, Topology  {\bf 25}  (1986),  no. 3, 375--380.

\bibitem[CFG2]{_CFG:_Frolicher_}
L.A. Cordero, M. Fern\'andez, A. Gray,  {\em The Fr\"olicher spectral sequence and complex compact nilmanifolds,} C. R. Acad. Sci. Paris
S\'er. I Math. {\bf 305} (1987), no. 17, 753--756.

\bibitem[CFL]{_CFL:_comple_}
L.A. Cordero, M. Fern\'andez, M. de Le\'on,  {\em Examples of
compact complex manifolds with no K\"ahler structure}, Portugal.
Math.  {\bf 44}  (1987),  no. 1, 49--62.



\bibitem[CFGU]{CFGU}
L.A. Cordero, M. Fern\'andez, A. Gray, L. Ugarte, {\em Compact
nilmanifolds with nilpotent complex structures: Dolbeault
cohomology}, Trans. Amer. Math. Soc. {\bf 352} (12), 5405--5433.

\bibitem[CG]{_Corwin_Greenlaf_}  L.J. Corwin, F.P. Greenleaf, {\em Representations of nilpotent Lie groups and their applications}, Cambridge studies in advanced mathematics 18, Cambridge, New York, 1990.


\bibitem[DP]{Demailly-Paun}
J.-P. Demailly, M. Paun, {\em Numerical characterization of the K\"ahler cone of a compact K\"ahler manifold}, Annals of Math. 159 (2004) 1247-1274.


\bibitem[DF1]{_Dotti_Fino:8-dim_}
I. Dotti, A. Fino,  {\em Abelian hypercomplex 8-dimensional
nilmanifolds,} Ann. Global Anal. Geom. {\bf 18} (1) (2000), 47--59.


\bibitem[FG]{_Fernandes_Gray:Iwasawa_}
M. Fern\'andez, A. Gray,
{\em The Iwasawa manifold},
Differential geometry, Pe\~n\'\i scola 1985, 157-159,
Lecture Notes in Math., 1209, Springer, Berlin, 1986.



\bibitem[FG]{_Fino_Gra_}
A. Fino,  G. Grantcharov,
{\em On some properties of the manifolds with
skew-symmetric torsion and holonomy $SU(n)$
and $Sp(n)$},
Adv. Math. 189 (2004), no. 2, 439--450.


%\bibitem[Fu83]{_Fujiki:structure_}
%Fujiki, A.,
%{\em On the structure of compact complex manifolds in C},
%Adv. Stud. Pure Math. 1, North-Holland, Amsterdam (1983) 231-302.

\bibitem[G]{_Guggenheimer:fail_}
H. Guggenheimer,
{\em Sur les vari\'et\'es qui poss\`edent une forme ext\'erieure quadratique
ferm\'ee,} C. R. Acad. Sci. Paris 232, (1951). 470-472.

 \bibitem[H]{_Hasegawa_}
K. Hasegawa,  {\em Minimal models of nilmanifolds, } Proc. Amer.
Math. Soc. {\bf 106} (1989), no. 1, 65--71.

\bibitem[Jac]{Jac} H. Jacobowitz,
{\em Transversely holomorphic foliations and CR structures.}
VI Workshop on Partial Differential Equations, Part I (Rio de Janeiro, 1999).
Mat. Contemp. {\bf 18} (2000), 175-194.


\bibitem[KS]{_Ketsetzis_Salamon_} G. Ketsetzis,  S.  Salamon, {\em Complex structures on the Iwasawa manifold}, Adv. Geom. 4 (2004), no. 2, 165-179.


\bibitem[M]{_Malcev_}
A. I. Mal'\v{c}ev, {\em On a class of homogeneous spaces}, AMS
Translation No. {\bf 39} (1951).

\bibitem[Mi]{Milnor}  J. Milnor, {\em  Curvatures  of  left  invariant  metrics  on Lie groups},  Adv.  Math.  {\bf 21} (1976),  293--329.

\bibitem[N]{Nomizu}  K. Nomizu, {\em On the cohomology of compact homogeneous spaces of nilpotent Lie groups}, Ann. of Math. {\bf 59} (1954), 531--538.

\bibitem[Pet]{Peternel} Th. Peternell,  {\em
  Modifications}, 
Several complex variables, VII, 285-317, Encyclopaedia Math. Sci., 74, Springer, Berlin, (1994)

\bibitem[Rag]{_Raghunatan_} M. S. Raghunathan, {\em Discrete subgroups of Lie groups}, Springer E.M.G. 68, Berlin, Heidelberg, New York, 1972.

\bibitem[R]{_Rollenske:def_at_large_}
S. Rollenske,
{\em Geometry of nilmanifolds with left-invariant complex
  structure and deformations in the large},
Proc. Lond. Math. Soc. (3) 99 (2009), no. 2, 425-460.

\bibitem[Sak]{_Sakane_} Y. Sakane, {\em On compact complex parallelalisable solvmanifolds}, Osaka J. Math. 13 (1976), 187-212.



\bibitem[S]{_Sal_}
S.M. Salamon,
{\em  Complex structures on nilpotent Lie algebras}, J. Pure Appl.
Algebra {\bf 157} (2001), 311--333.


\bibitem[Sam]{Samelson} H. Samelson, {\em A class of compact-analytic manifolds}, Portugaliae Math. 12 (1953),  129-132.


\bibitem[T]{_Thurston:Kodaira_}
W. P. Thurston,  {\em Some simple examples of symplectic manifolds,}
Proc. Amer. Math. Soc. 55 (1976), no. 2, 467-468.

\bibitem[Ue75]{_Ueno_}
 K. Ueno,
{\em Classification theory of algebraic varieties and
compact complex spaces,}
Lecture Notes in Mathematics 439, Springer-Verlag, Berlin-New York (1975).

\bibitem[Ug]{_Ugarte_} L. Ugarte,
{\em Hermitian structures on six-dimensional nilmanifolds}, Transform.  Groups {\bf 12} (2007), 175--202.

\bibitem[W]{_Winkelmann:paralleli_}
J. Winkelmann,
{\em  Complex Analytic Geometry of Complex Parallelizable Manifolds,}
Mem. Soc. Math. France (N.S.) 72-73, Soc. Math. France, Montrouge, 1998.

\end{thebibliography}
}
\end{document}